\documentclass[final,3p,times]{elsarticle}
\usepackage{amsmath} 
\usepackage{amssymb}
\usepackage{textcomp}
\usepackage{graphicx}		
\usepackage{float}
\usepackage{color}
\usepackage{url}
\usepackage{array}
\usepackage{booktabs,ctable,multirow}
\usepackage{bm}
\usepackage[T1]{fontenc}
\usepackage{epsfig}
\usepackage{epstopdf}
\usepackage{url}
\usepackage{color}
\newcommand{\red}{\color[rgb]{1,0,0}}
\usepackage{caption}
\usepackage{subcaption}

\newcommand{\R}{\mathbb{R}}
\newcommand{\xx}{ {\bm x} }
\newcommand{\yy}{ {\bm y} }
\newcommand{\zz}{ {\bm z} }
\newcommand{\ww}{ {\bm w} }
\newcommand{\kk}{ {\bm k} } 
 
\newcommand{\bfone}{\bm{1}}
\newcommand{\bfi}{\bm{\phi}}
\newcommand{\bFi}{\bm{\Phi}}

\newcommand{\hloc}{{\overline{h}_{\text{local}}}}
\newtheorem{definition}{Definition}

\begin{document}

\begin{frontmatter}
  \title{Inverting Nonlinear Dimensionality Reduction with\\
    Scale-Free Radial Basis Function Interpolation}
  \author[add1]{Nathan D. Monnig\corref{cor1}}
  \ead{E-mail: nathan.monnig@colorado.edu}
  \ead[url]{http://amath.colorado.edu/student/monnign/}
  \author[add1]{Bengt Fornberg}
  \author[add2]{Fran\c{c}ois G. Meyer}
  \cortext[cor1]{Corresponding author}
  \address[add1]{Department of Applied Mathematics, UCB 526, University of
    Colorado at Boulder, Boulder, CO  80309}
  \address[add2]{Department of Electrical Engineering, UCB 425, University of
    Colorado at Boulder, Boulder, CO  80309}
  \begin{abstract}
    Nonlinear dimensionality reduction embeddings computed from datasets do not
    provide a mechanism to compute the inverse map. In this paper, we
    address the problem of computing a stable inverse map to such a
    general bi-Lipschitz map. Our approach relies on radial basis
    functions (RBFs) to interpolate the inverse map everywhere on the
    low-dimensional image of the forward map. We demonstrate that the
    scale-free cubic RBF kernel performs better than the Gaussian
    kernel: it does not suffer from ill-conditioning, and does not
    require the choice of a scale.  The proposed construction is shown
    to be similar to the Nystr{\"o}m extension of the eigenvectors of
    the symmetric normalized graph Laplacian matrix.  Based on this
    observation, we provide a new interpretation of the Nystr{\"o}m
    extension  with suggestions for improvement.
  \end{abstract}
\end{frontmatter}

\section{Introduction}
\label{intro}
The construction of parametrizations of low dimensional data in high
dimension is an area of intense research (e.g., 
\cite{belkin2003,Coifman06b,roweis,scholkopf}).  A major
limitation of these methods is that they are only defined on a
discrete set of data.  As a result, the inverse mapping is also only
defined on the data.  There are well known strategies to extend the
forward map to new points---for example, the Nystr{\"o}m extension is a
common approach to solve this {\em out-of-sample extension} problem
(see e.g., \cite{coifman} and references therein).  However, the
problem of extending the inverse map (i.e. the {\em preimage problem})
has received little attention so far (but see \cite{elgammal2008}).
The nature of the preimage problem precludes application of the
Nystr{\"o}m extension, since it does not involve extension of
eigenvectors.

We present a method to numerically invert a general smooth
bi-Lipschitz nonlinear dimensionality reduction mapping over all
points in the image of the forward map.  The method relies on
interpolation via radial basis functions of the coordinate functions
that parametrize the manifold in high dimension.

The contributions of this paper are twofold.  Primarily, this paper
addresses a fundamental problem for the analysis of datasets: given
the construction of an adaptive parametrization of the data in terms
of a small number of coordinates, how does one synthesize new data
using new values of the coordinates? We provide a simple and elegant
solution to solve the ``preimage problem''. Our approach is scale-free
and numerically stable and can be applied to any nonlinear dimension
reduction technique. The second contribution is a novel interpretation
of the Nystr{\"o}m extension as a properly rescaled radial basis function interpolant.
A precise analysis of this similarity yields a critique of the
Nystr{\"o}m extension, as well as suggestions for improvement.

\section{The Inverse Mapping}
\subsection{Definition of the problem, and approach}
\label{inverse}
We consider a finite set of $n$ datapoints $\{ \xx^{(1)}, \ldots ,
\xx^{(n)} \} \subset \R^D$ that lie on a bounded low-dimensional
smooth manifold ${\cal M} \subset \R^D$, and we assume that a
nonlinear mapping has been defined for each point $\xx^{(i)}$,
\begin{align}
  \bFi_n: {\cal M} \subset \R^D & \longrightarrow \R^d\\
  \xx^{(i)} & \longmapsto \yy^{(i)} = \bFi_n(\xx^{(i)}), \quad  i=1,\ldots,n. 
\end{align}
We further assume that the map $\bFi_n$ converges toward a limiting
continuous function, $\bFi: {\cal M} \rightarrow \bFi\left({\cal
    M}\right)$, when the number of samples goes to infinity. Such
limiting maps exist for algorithms such as the Laplacian eigenmaps
\cite{Coifman06b}.

In practice, the construction of the map $\bFi_n$ is usually only the
first step. Indeed, one is often interested in exploring the
configuration space in $\R^d$, and one needs an inverse map to
synthesize a new measurement $\xx$ for a new configuration $\yy
= \begin{bmatrix}y_1 & \cdots & y_d \end{bmatrix}^T$ in the coordinate
domain (see e.g., \cite{Coifman08}).  In other words, we would like to
define an inverse map $\bFi_n^{-1} (\yy)$ at any point $\yy \in
\bFi_n\left({\cal M}\right)$.  Unfortunately, unlike linear methods
(such as PCA), nonlinear dimension reduction algorithms only provide
an explicit mapping for the original discrete dataset $\{ \xx^{(1)},
\ldots , \xx^{(n)} \}$. Therefore, the inverse mapping $\bFi_n^{-1}$
is only defined on these data.

The goal of the present work is to generate a numerical extension of
$\bFi_n^{-1}$ to all of $\bFi({\cal M}) \subset \R^d$.  To simplify
the problem, we assume the mapping $\bFi_n$ coincides with the
limiting map $\bFi$ on the data, $\bFi_n(\xx^{(i)}) = \bFi(\xx^{(i)})$
for $i=1,\ldots,n$.  This assumption allows us to rephrase the problem
as follows: we seek an extension of the map $\bFi^{-1}$ everywhere on
$\bFi({\cal M})$, given the knowledge that $\bFi^{-1}(\yy^{(i)})
=\bFi_n^{-1}(\yy^{(i)}) = \xx^{(i)}$.  We address this problem using
interpolation, and we construct an approximate inverse
$\bFi_n^\dagger$, which converges toward the true inverse as the
number of samples, $n$, goes to infinity,
\begin{align}
  & \bFi^\dagger: \bFi\left({\cal M}\right) \rightarrow \R^D,\quad
  \text{with} \quad
  \bFi^\dagger \left(\yy^{(i)}\right)  = \xx^{(i)}, \\
  \text{and} & \quad \forall \yy \in \bFi\left({\cal M}\right), \quad
  \lim_{n \rightarrow \infty} \bFi_n^\dagger (\yy) = \bFi^{-1}(\yy).
\end{align}
Using terminology from geometry, we call $\bFi({\cal M})$ the {\em
  coordinate domain}, and $\bFi^{-1}$ a {\em coordinate map} that
parametrizes the manifold  ${\cal M} = \{\xx \in \R^D; \xx = \bFi^{-1}(\yy) ,
\yy \in \bFi ({\cal M})\}$. The components of $\bFi^{-1}
= \begin{bmatrix} \phi_1^{-1} \cdots \phi_D^{-1} \end{bmatrix}^T$
are the {\em coordinate functions}.  We note that the
focus of the paper is not the construction of new points $\yy$ in the
coordinate domain, but rather the computation of the coordinate
functions everywhere in $\bFi ({\cal M})$.

\subsection{Interpolation of multivariate functions defined on scattered data}
\label{rbf}
Given the knowledge of the inverse at the points $\yy^{(i)}$, we wish
to interpolate $\bFi^{-1}$ over $\bFi\left({\cal M}\right)$. We
propose to interpolate each coordinate function, $\phi_i^{-1} (\yy), i
=1,\ldots, D$ independently of each other. We are thus facing the
problem of interpolating a function of several variables defined on
the manifold $\bFi({\cal M})$. Most interpolation techniques that are
designed for single variable functions can only be extended using
tensor products, and have very poor performance in several
dimensions. For instance, we know from Mairhuber theorem (e.g.,
\cite{fasshauer}) that we should not use a basis independent of the
nodes (for example, polynomial) to interpolate scattered data in
dimension $d>1$.  As a result, few options exist for multivariate
interpolation. Some of the most successful interpolation methods
involve Radial Basis Functions (RBFs) \cite{fasshauer}. Therefore, we
propose to use RBFs to construct the inverse mapping. Similar methods
have been explored in \cite{elgammal2008,powell} to interpolate data
on a low-dimensional manifold. We note that while kriging \cite{wackernagel} is
another
common approach for interpolating scattered data, most kriging
techniques are equivalent to RBF interpolants \cite{scheuerer}. In fact,
because in our application we lack specialized information about the
covariance structure of the inverse map, kriging is identical to RBF
interpolation.

We focus our attention on two basis functions: the Gaussian and the
cubic.  These functions are representative of the two main classes of
radial functions: scale dependent, and scale invariant.  In the
experimental section we compare the RBF methods to Shepard's method
\cite{shepard}, an approach for multivariate interpolation and
approximation that is used extensively in computer graphics
\cite{lewis}, and which was recently proposed in \cite{kushnir} to
compute a similar inverse map.

For each coordinate function $\phi_i^{-1}$, we define $\phi^\dagger_i$
to be the RBF interpolant to the data $\left(\yy^{(i)}, \xx^{(i)}\right)$,
\begin{equation}
  \text{for all}\;  \yy \in \bFi({\cal M}), \quad \phi^\dagger_i
  (\yy)=\sum_{j=1}^n \alpha_i^{(j)} k(\yy,\yy^{(j)}).
  \label{one-d}
\end{equation}
The reader will notice that we dropped the dependency on $n$ (number
of samples) in $\bFi^\dagger= \begin{bmatrix} \phi^\dagger_1 \ldots
  \phi^\dagger_D \end{bmatrix}^T$ to ease readability. The function $k$
in (\ref{one-d}) is the kernel that defines the radial basis
functions, $k(\zz, \ww) = g(\|\zz - \ww\|)$. The weights, $\{
\alpha_i^{(1)},\ldots,\alpha_i^{(n)} \}$, are determined by imposing
the fact that the interpolant be exact at the nodes
$\yy^{(1)},\ldots,\yy^{(n)}$, and thus are given by the solution of
the linear system
\begin{equation}
  \begin{bmatrix}
    k(\yy^{(1)},\yy^{(1)} ) &
    \cdots & k( \yy^{(1)},\yy^{(n)} ) \\
    \vdots & \ddots & \vdots \\
    k(\yy^{(n)},\yy^{(1)} ) & \cdots  & k(\yy^{(n)},\yy^{(n)} )
  \end{bmatrix}
  \begin{bmatrix}
    \alpha_i^{(1)} \\
    \vdots \\
    \alpha_i^{(n)} 
  \end{bmatrix}
  =    
  \begin{bmatrix}
    x_i^{(1)} \\
    \vdots \\
    x_i^{(n)}
  \end{bmatrix}.
  \label{one-coord}
\end{equation}
We can combine the $D$ linear systems (\ref{one-coord}) by concatenating all the
coordinates in the right-hand side of (\ref{one-coord}), and the
corresponding unknown weights on the left-hand side of
(\ref{one-coord}) to form the system of equations,
\begin{equation}
  \begin{bmatrix}
    k( \yy^{(1)},\yy^{(1)} ) & \cdots & k( \yy^{(1)},\yy^{(n)} ) \\ 
    \vdots & \ddots & \vdots \\
    k( \yy^{(n)},\yy^{(1)} ) & \cdots & k( \yy^{(n)},\yy^{(n)} )
  \end{bmatrix}
  \begin{bmatrix}
    \alpha_1^{(1)} &  & \alpha_D^{(1)} \\
    \vdots & \cdots & \vdots \\
    \alpha_1^{(n)} &  & \alpha_D^{(n)} 
  \end{bmatrix}
  =
  \begin{bmatrix}
    x_1^{(1)} &  & x_D^{(1)} \\
    \vdots & \cdots & \vdots \\
    x_1^{(n)} &  & x_D^{(n)}
  \end{bmatrix},
  \label{gauss_rbf}
\end{equation}
which takes the form $KA = X$, where $K_{i,j} = k (\yy^{(i)},
\yy^{(j)})$, $A_{i,j} = \alpha_j^{(i)}$, and $X_{i,j} =
x_j^{(i)}$. Let us define the vector $\kk(\yy,\cdot)=\begin{bmatrix}
  k(\yy,\yy^{(1)}) & \ldots &k(\yy,\yy^{(n)})
\end{bmatrix}^T$. The approximate inverse at a point $\yy \in
\bFi({\cal M})$ is given by
\begin{equation}
  \bFi^\dagger(\yy)^T = \kk(\yy,\cdot)^T A =
  \kk(\yy,\cdot)^T K^{-1} X.
  \label{rbf_eqn}
\end{equation}
\section{Convergence of RBF Interpolants}
\subsection{Invertibility and Conditioning}
\label{conditioning}
The approximate inverse (\ref{rbf_eqn}) is obtained by interpolating
the original data $\left(\yy^{(i)}, \xx^{(i)}\right)$ using RBFs. In
order to assess the quality of this inverse, three questions must be
addressed: 1) Given the set of interpolation nodes,
$\left\{\yy^{(i)}\right\}$, is the interpolation matrix $K$ in
(\ref{gauss_rbf}) necessarily non-singular and well-conditioned?  2)
How well does the interpolant (\ref{rbf_eqn}) approximate the true inverse
$\bFi^{-1}$? 3) What convergence rate can we expect as we populate
the domain with additional nodes?  In this section we provide elements
of answers to these three questions.  For a detailed treatment, see
\cite{fasshauer,wendland}.

In order to interpolate with a radial basis function $k(\zz,\ww) =
g(\|\zz - \ww\|)$, the system (\ref{gauss_rbf}) should have a unique
solution and be well-conditioned. In the case of the Gaussian defined
by
\begin{equation}
  k(\zz,\ww) = \exp(-\varepsilon^2 \| \zz - \ww \|^2),
\end{equation}
the eigenvalues of $K$ in (\ref{gauss_rbf}) follow patterns in the
powers of $\varepsilon$ that increase with successive eigenvalues,
which leads to rapid ill-conditioning of $K$ with increasing $n$
(e.g., \cite{fornberg2007runge}; see also \cite{bermanis} for a
discussion of the numerical rank of the Gaussian kernel). The
resulting interpolant will exhibit numerical {\it saturation error}.
This issue is common among many scale-dependent RBF interpolants.  The
Gaussian scale parameter, $\varepsilon$, must be selected to match the spacing of the
interpolation nodes.  One commonly used measure of node spacing is the
{\em fill distance}, the maximum distance from an interpolation node.
\noindent 
\begin{definition}
  For the domain $\Omega \subset \R^d$ and a
  set of interpolation nodes
  $Z = \{\zz^{(1)},\ldots,\zz^{(n)} \} \subset \Omega$ the {\em fill
    distance}, $h_{Z,\Omega}$, is defined by
  \begin{equation}
    h_{Z,\Omega} := \sup_{\zz \in \Omega} \min_{\zz^{(j)} \in Z} \lVert \zz -
    \zz^{(j)} \rVert.
  \end{equation}
\end{definition}
\begin{figure}[H]
  \begin{center}
    \includegraphics[width=0.8\textwidth]{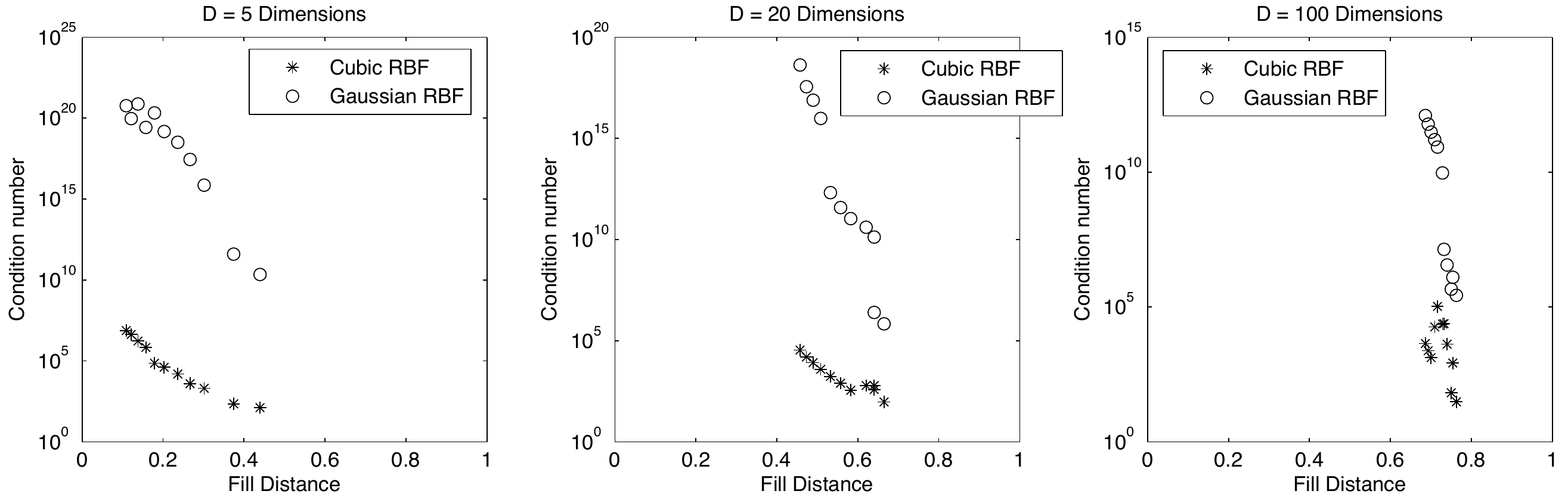}
  \end{center}
  \caption{Condition number of $K$ in (\ref{gauss_rbf}), for the Gaussian
    ($\circ$) and the cubic ($\ast$) as a function of the fill distance
    for a fixed scale $\varepsilon = 10^{-2}$.  Points are randomly
    scattered on the first quadrant of the unit sphere in $\R^D$, for
    $D =5, 20, 100$ from left to right.  Note: the same range of $n$,
    from 10 to 1000, was used in each dimension.  In high dimension,
    it takes a large number of points to reduce fill distance.
    However, the condition number of $K$ still grows rapidly for
    increasing $n$.  }
  \label{cond_W_fill}
\end{figure}
\begin{figure}[H]
  \begin{center}
    \includegraphics[width=0.8\textwidth]{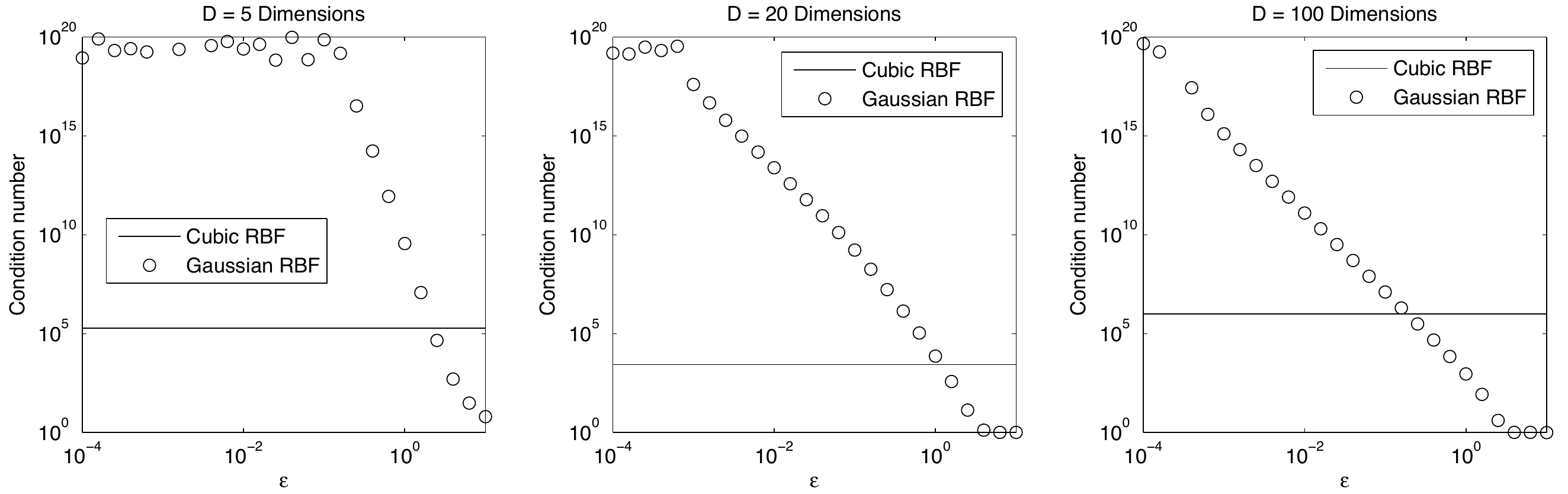}
  \end{center}
  \caption{Condition number of $K$ in (\ref{gauss_rbf}), for the Gaussian ($\circ$)
    and the cubic (--) as a function of the scale $\varepsilon$, for a fixed
    fill distance. $n=200$ points are randomly scattered on the first
    quadrant of the unit sphere in $\R^D$, $D =5, 20, 100$ from left to right.
    \label{cond_W_ep}}
\end{figure}
\noindent Owing to the difficulty in precisely establishing the boundary of a
domain $\Omega \subset \R^d$ defined by a discrete set of sampled
data, estimating the fill distance $h_{Z,\Omega}$ is somewhat
difficult in practice.  Additionally, the fill distance is a measure
of the ``worst case'', and may not be representative of the
``typical'' spacing between nodes.  Thus, we consider a proxy for fill
distance which depends only on mutual distances between the data
points.  We define the  {\it local fill distance}, $\hloc$, to denote
the average distance to a nearest neighbor,
\begin{equation}
  \hloc := \frac{1}{n} \sum_{i=1}^n \min_{j \neq i} \lVert \zz^{(i)} -
  \zz^{(j)}
  \rVert.
  \label{localfill}
\end{equation}
The relationship between the condition number of $K$ and the spacing
of interpolation nodes is explored in Fig.~\ref{cond_W_fill}, where we
observe rapid ill-conditioning of $K$ with respect to decreasing local
fill distance, $\hloc$. Conversely, if $\hloc$ remains constant while
$\varepsilon$ is reduced, the resulting interpolant improves until
ill-conditioning of the $K$ matrix leads to propagation of numerical
errors, as is shown in Fig.~\ref{cond_W_ep}. When interpolating with
the Gaussian kernel, the choice of the scale parameter $\varepsilon$
is difficult. On the one hand, smaller values of $\varepsilon$ likely
lead to a better interpolant. For example, in 1-$d$, a Gaussian RBF
interpolant will converge to the Lagrange interpolating polynomial in
the limit as $\varepsilon \rightarrow 0$ \cite{fornberg2002limit}. On
the other hand, the interpolation matrix becomes rapidly
ill-conditioned for decreasing $\varepsilon$. While some stable
algorithms have been recently proposed to generate RBF interpolants
(e.g., \cite{fornberg2013}, and references therein) these
sophisticated algorithms are more computationally intensive and
algorithmically complex than the RBF-Direct method used in this paper,
making them undesirable for the inverse-mapping interpolation task.

Saturation error can be avoided by using the scale-free RBF kernel
$g(\|\zz -\ww\|)= \lVert \zz-\ww \rVert^3$, one
instance from the set of RBF kernels known as the {\it radial powers},
\begin{equation}
  \begin{array}{ccc}
    g(\|\zz-\ww\|)= \lVert \zz-\ww \rVert^\rho & \text{for} & \rho = 1,3,5,\ldots. 
  \end{array}
  \label{radial-powers}
\end{equation}
\noindent Together with the {\it thin plate splines},
\begin{equation}
  \begin{array}{ccc}
    g(\|\zz - \ww\|)= \lVert \zz-\ww \rVert^\rho \log \lVert \zz-\ww \rVert & \text{for}
    &
    \rho = 2,4,6,\ldots,
  \end{array}
\end{equation}
\noindent \noindent they form the family of RBFs known as the {\it
  polyharmonic splines}.  

Because it is a monotonically increasing function, the cubic kernel,
$\lVert \zz-\ww \rVert^3$, may appear less intuitive than the
Gaussian. The importance of the cubic kernel stems from the fact that
the space generated by linear combinations of shifted copies of
the kernel is composed of splines. In one dimension, one recovers the
cubic spline interpolant.  One should note that the behavior of the
interpolant in the far field (away from the boundaries of the convex
hull of the samples) can be made linear (by adding constants and
linear polynomials) as a function of the distance, and therefore
diverges much more slowly than $r^3$ \cite{fornberg2002boundaries}.

In order to prove the existence and uniqueness of an interpolant of
the form, 
\begin{equation}
  \phi^\dagger_i(\yy)=\sum_{j=1}^n \alpha_i^{(j)} \| \yy - \yy^{(j)} \|^3 + \gamma_i +
  \sum_{k=1}^{d} \beta_{k,i} y_k,
\end{equation}
we require that the set $\{\yy^{(1)},\ldots,\yy^{(n)} \}$ be a
$1$-unisolvent set in $\R^d$, where {\em $m$-unisolvency} is as
follows.
\begin{definition}
  The set of nodes $\{\zz^{(1)},\ldots,\zz^{(n)} \} \subset \R^d$ is
  called {\em $m$-unisolvent} if the unique polynomial of total degree
  at most $m$ interpolating zero data on $\{ \zz^{(1)},\ldots,
  \zz^{(n)} \}$ is the zero polynomial.
\end{definition}
For our problem, the condition that the set of nodes $\{\yy^{(j)}\}$
be 1-unisolvent is equivalent to the condition that the matrix 
\begin{equation}
  \begin{bmatrix}
    1 & \cdots & 1\\
    \begin{bmatrix}
      \\
      \yy^{(1)}
      \\
      \\
    \end{bmatrix}
     & \cdots & 
    \begin{bmatrix}
      \\
      \yy^{(n)}
      \\
      \\
    \end{bmatrix}
  \end{bmatrix}
  =
  \begin{bmatrix}
    1 & \cdots & 1\\
    \phi_1 (\xx_1) & \cdots & \phi_1 (\xx_n)\\
    \vdots & & \vdots\\
    \phi_d (\xx_1) & \cdots & \phi_d (\xx_n)\\
  \end{bmatrix}
  \label{polymatrix}
\end{equation}
have rank $d+1$ (we assume that $n \geq d+1$).  This condition is
easily satisfied. Indeed, the rows $2,\ldots,d+1$ of
(\ref{polymatrix}) are formed by the orthogonal eigenvectors of
$D^{-1/2} W D^{-1/2}$.  Additionally, the first eigenvector, $\bfi_0$,
has constant sign. As a result, $\bfi_1,\ldots,\bfi_d$ are linearly
independent of any other vector of constant sign, in particular
$\bfone$.  In Figures \ref{cond_W_fill} and \ref{cond_W_ep} we see that
the cubic RBF system exhibits much better conditioning than the
Gaussian.

\subsection{What Functions Can We Reproduce? the Concept of Native Spaces}
We now consider the second question: can the interpolant (\ref{one-d})
approximate the true inverse $\bFi^{-1}$ to arbitrary precision?  As
we might expect, an RBF interpolant will converge to functions
contained in the completion of the space of linear combinations of the
kernel, ${\cal H}_k(\Omega) = span\{k(\cdot,\zz) : \zz \in \Omega \}$.
This space is called the {\em native space}. We note that the
completion is defined with respect to the $k$-norm, which is induced
by the inner-product given by the reproducing kernel $k$ on the
pre-Hilbert space ${\cal H}_k(\Omega)$ \cite{fasshauer}.

It turns out that the native space for the Gaussian RBF is a very
small space of functions whose Fourier transforms decay faster than a
Gaussian \cite{fasshauer}. In practice, numerical issues usually
prevent convergence of Gaussian RBF interpolants, even within the
native space, and therefore we are not concerned with this issue. The
native space of the cubic RBF, on the other hand, is an extremely
large space. When the dimension, $d$, is odd, the native space of the
cubic RBF is the Beppo Levi space on $\R^d$ of order $l=(d+3)/2$
\cite{wendland}. We recall the definition of a {\em Beppo Levi
  space} of order $l$.
\begin{definition}
  For $l>d/2$, the linear space
  $BL_l(\R^d) := \{f \in C(\R^d) : D^\alpha f \in  L_2(\R^d), \forall |\alpha|=l \}$,
  equipped with the inner product
  $\langle f,g \rangle_{\text{BL}_l(\R^d)} = \sum_{|\alpha|=l}
  \frac{l!}{\alpha!} \langle D^{\alpha} f, D^{\alpha} g
  \rangle_{L_2(\R^d)}$,
  is called the {\em Beppo Levi space on $\R^d$ of order $l$}, where
  $D^{\alpha}$ denotes the weak derivative of (multi-index) order $\alpha \in
  \mathbb{N}^d$ on $\R^d$.
\end{definition}
For even dimension, the Beppo Levi space on $\R^d$ of order
$l=(d+2)/2$ corresponds to the native space of the thin plate spline
$g(\|\zz - \ww\|)= \lVert \zz-\ww \rVert^2 \log \lVert \zz-\ww \rVert$
\cite{wendland}. Because we assume that the inverse map $\bFi^{-1}$ is
smooth, we expect that it belongs to any of the Beppo Levi
spaces. Despite the fact that we lack a theoretical characterization
of the native space for the cubic RBF in even dimension, all of our
numerical experiments have demonstrated equal or better performance of
the cubic RBF relative to the thin plate spline in all dimensions (see
also \cite{wild} for similar conclusions).  Thus, to promote
algorithmic simplicity for practical applications, we have chosen to
work solely with the cubic RBF.
\subsection{Convergence Rates}
The Gaussian RBF interpolant converges (in $L^\infty$ norm)
exponentially fast toward functions in the native space, as a function
of the decreasing fill distance $h_{Z,\Omega}$ \cite{fasshauer}.
However, as observed above, rapid ill-conditioning of the
interpolation matrix makes such theoretical results irrelevant without
resorting to more costly stable algorithms.  The cubic interpolant
converges at least as fast as ${\cal O}(h_{Z,\Omega}^{3/2})$ in the
respective native space \cite{wendland}. In practice, we have
experienced faster rates of algebraic convergence, as shown in the
experimental section.
\section{Experiments}
\label{experiments}
We first conduct experiments on a synthetic manifold, and we then
provide evidence of the performance of our approach on real data.  For
all experiments we quantify the performance of the interpolation using
a ``leave-one-out reconstruction'' approach: we compute
$\bFi^{\dagger}\left(\yy^{(j)}\right)$, for $j=1,\ldots,n$, using the
remaining $n-1$ points: $\{\yy^{(1)}, \ldots,\yy^{(j-1)}$, $\yy^{(j+1)},\ldots,\yy^{(n)} \}$ and
their coordinates in $\R^D$, $\{
\xx^{(1)},\ldots,\xx^{(j-1)},\xx^{(j+1)},\ldots,\xx^{(n)} \}$. The
average performance is then measured using the average leave-one-out
$l^2$ reconstruction error,
\begin{equation}
  E_{\text{avg}} = \frac{1}{n}\sum_{j=1}^n \lVert \xx^{(j)} -
  \bFi^{\dagger}(\yy^{(j)}) \rVert.
  \label{error}
\end{equation}
In order to quantify the effect of the sampling density on the
reconstruction error, we compute $E_\text{avg}$ as a function of
$\hloc$, which is defined by (\ref{localfill}). The two RBF
interpolants are compared to Shepard's method, a multivariate
interpolation/approximation method used extensively in computer
graphics \cite{lewis}. Shepard's method computes the optimal constant
function that minimizes the sum of squared errors within a
neighborhood ${\cal N}_\yy$ of $\yy$ in $\R^d$, weighted according to
their proximity to $\yy$. The solution to this moving least squares
approximation is given by
\begin{equation}
  \bFi^\dagger_{\text{Shepard}} (\yy) = \sum_{j:\yy^{(j)} \in {\cal
      N}_\yy} \frac{\exp(-\varepsilon^2 \lVert \yy-\yy^{(j)} \rVert^2
    ) }{ \sum_{i:\yy^{(i)} \in {\cal N}_\yy}  \exp(-\varepsilon^2
    \lVert \yy-\yy^{(i)} \rVert^2 )} \xx^{(j)}.
\end{equation}
The relative impact of neighboring function values is controlled by
the scale parameter $\varepsilon$, which we choose to be a multiple of
$1/\hloc$.
\subsection{Unit Sphere in $\R^D$}
\begin{figure}[H]
  \begin{center}
    \includegraphics[width=0.8\textwidth]{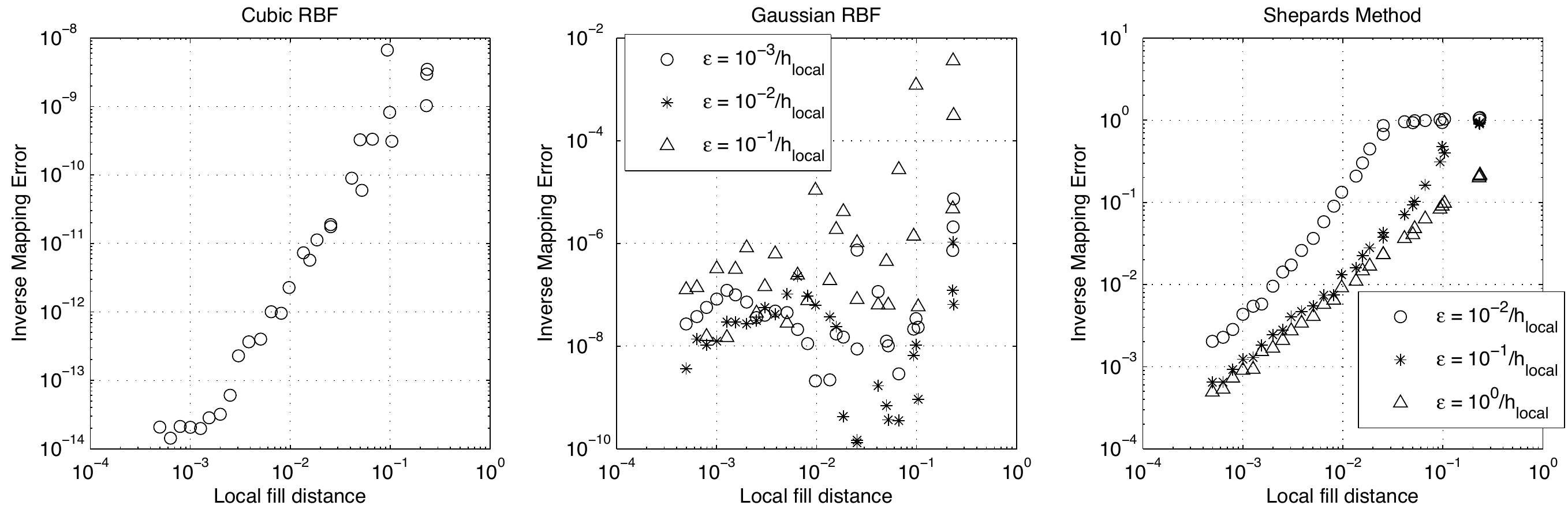}
  \end{center}
  \caption{Average leave-one-out reconstruction residual, $E_\text{avg}$, on $S^4$
    embedded in $\R^{10}$, using the cubic (left), the Gaussian
    (center), and Shepard's method (right).  Note the difference in
    the range
    of $y$-axis.}
  \label{unitsphere}
\end{figure}

For our synthetic manifold example, we sampled $n$ points $\{
\xx^{(1)},\ldots,\xx^{(n)} \}$ from the uniform distribution on the unit sphere
$S^4$, then embedded these data in $\R^{10}$ via a random unitary
transformation.
The data are mapped to $\{ \yy^{(1)},\ldots,\yy^{(n)} \} \subset
\R^d=\R^5$ using the first five non-trivial eigenvectors of the graph
Laplacian.  The minimum of the total number of available neighbors, $n-1$, and 200
neighbors was used to compute the interpolant.  For each local fill
distance, $\hloc$, the average reconstruction error is computed using
(\ref{error}).  The performances of the cubic RBF, Gaussian RBF, and
Shepard's method versus $\hloc$ are shown in Fig.~\ref{unitsphere}. We
note that the interpolation error based on the cubic RBF is lowest,
and appears to scale approximately with $O(\overline{h}^2_{\text{local}})$, an improvement
over the ${\cal O}(h_{Z,\Omega}^{3/2})$ bound \cite{wendland}.  In
fact, the cubic RBF proves to be extremely accurate, even with
a very sparsely populated domain: the largest $\hloc$ corresponds to
10 points scattered on $S^4$.
\noindent 
\begin{table}[H]
  \begin{center}
    \begin{small}
      \begin{tabular}{cccccccccccc}
        \toprule
        & scale ($\varepsilon \times \hloc$)  & {\bf 0} & {\bf 1} &
        {\bf 2} & {\bf 3} & {\bf 4} & {\bf 5} & {\bf 6} & {\bf 7} & {\bf 8} &
        {\bf 9}\\
        \midrule
        Cubic & ---                    & \red 0.248 & \red 0.135 & \red 0.349
        & \red 0.334 & \red 0.299 & \red 0.350 & \red 0.259 & \red 0.261 &
        \red 0.354 & \red 0.262 \\ 
        \midrule
        & $0.5$ & 0.319 & 0.169 & 0.421 & 0.417 & 0.362 & 0.424 & 0.315 & 0.314 & 0.452 & 0.313 \\
        Gaussian & $1$   & 0.305 & 0.223 & 0.375 & 0.363 & 0.345 & 0.382 & 0.310 & 0.322 & 0.369 & 0.312 \\
        & $2$   & 0.457 & 0.420 & 0.535 & 0.505 & 0.513 & 0.554 & 0.477 &
        0.497 & 0.491 & 0.478 \\
        \midrule
        & $0.5$ & 0.422 & 0.271 & 0.511 & 0.489 & 0.475 & 0.508 & 0.439 & 0.453 & 0.474 & 0.434 \\
        Shepard  & $1$   & 0.302 & 0.175 & 0.396 & 0.385 & 0.348 & 0.378 & 0.314 & 0.318 & 0.379 & 0.309 \\
        & $2$   & 0.303 & 0.186 & 0.400 & 0.382 & 0.362 & 0.402 & 0.320 & 0.325 & 0.382 & 0.320 \\
        \bottomrule
      \end{tabular}
    \end{small}
  \end{center}
  \caption{Reconstruction error $E_{\text{avg}}$ for each
       digit  (0-9).  Red denotes lowest average
    reconstruction residual.}
  \label{digits_table}
\end{table}

\begin{figure}	
	\centering
	\begin{subfigure}[t]{0.89in}
		\centering
\caption*{\scriptsize Original}\vspace{-0.09in}
		\includegraphics[width=0.89in,
trim=0.5in 0.4in 0.35in 0.3in,clip=true]{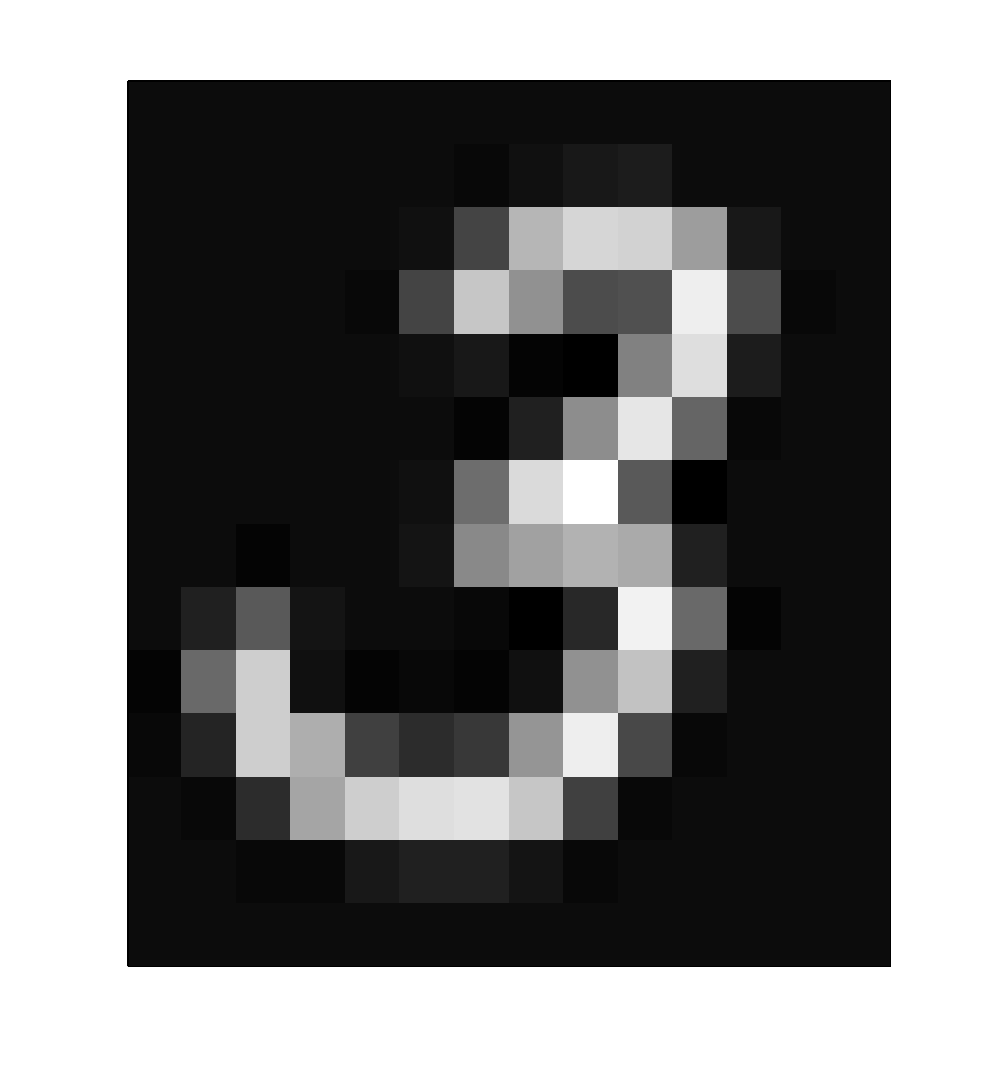}
	\end{subfigure}
	\begin{subfigure}[t]{0.89in}
		\centering
\caption*{\scriptsize Cubic RBF}\vspace{-0.09in}
		\includegraphics[width=0.89in,
trim=0.5in 0.4in 0.35in 0.3in,clip=true]{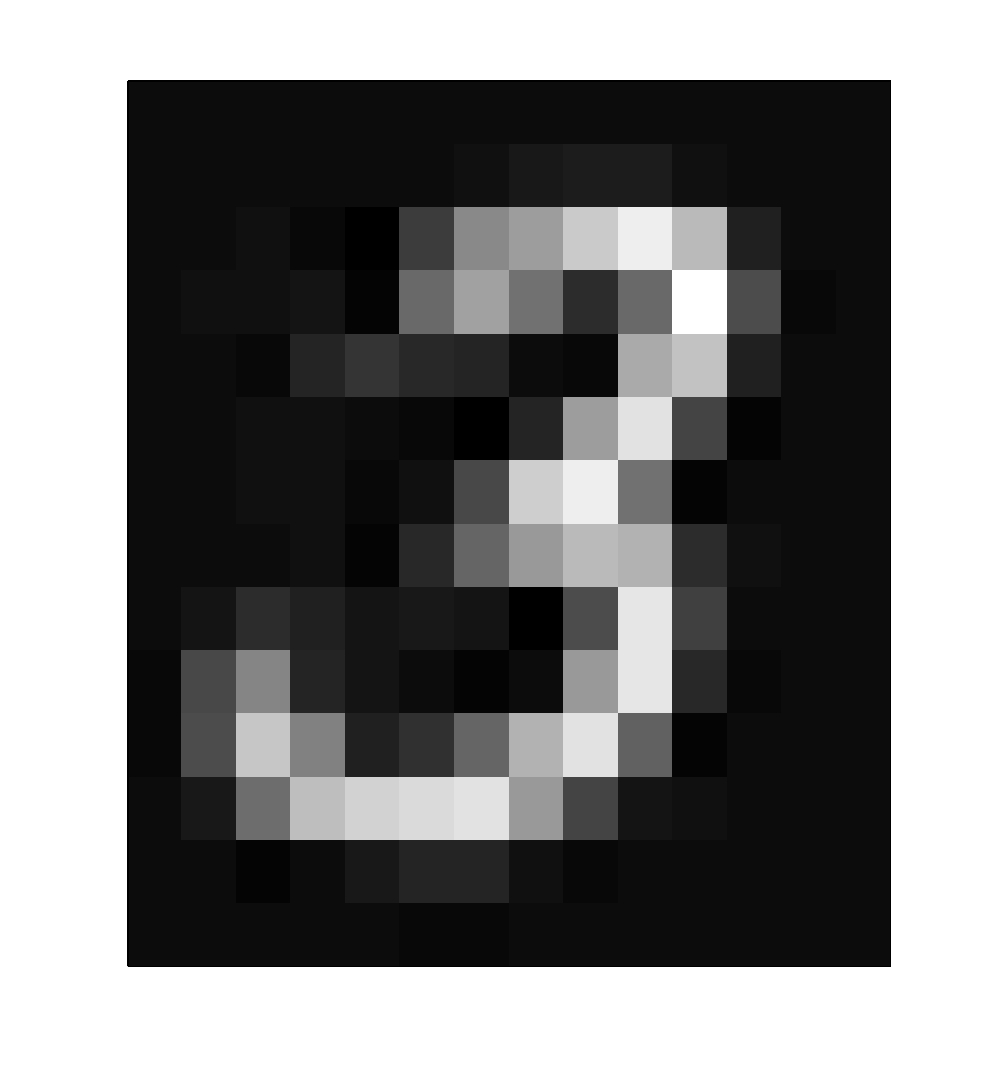}
	\end{subfigure}
	\begin{subfigure}[t]{0.89in}
		\centering
\caption*{\scriptsize Error = 0.23}\vspace{-0.09in}
		\includegraphics[width=0.89in,
trim=0.5in 0.4in 0.35in 0.3in,clip=true]{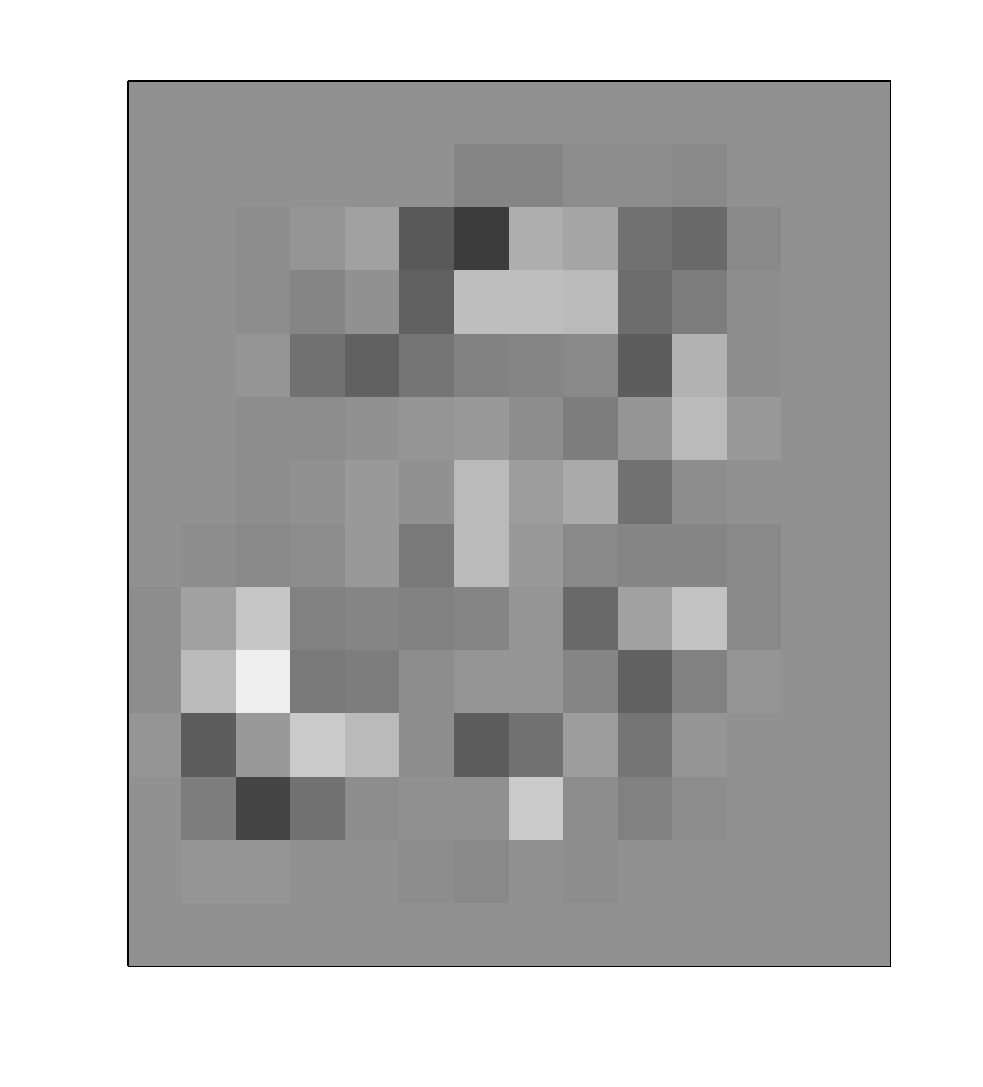}
	\end{subfigure}
	\begin{subfigure}[t]{0.89in}
		\centering
\caption*{\scriptsize Gaussian \tiny ($\varepsilon = 1/\bar
h_{local}$)}\vspace{-0.082in}
		\includegraphics[width=0.89in,
trim=0.5in 0.4in 0.35in 0.3in,clip=true]{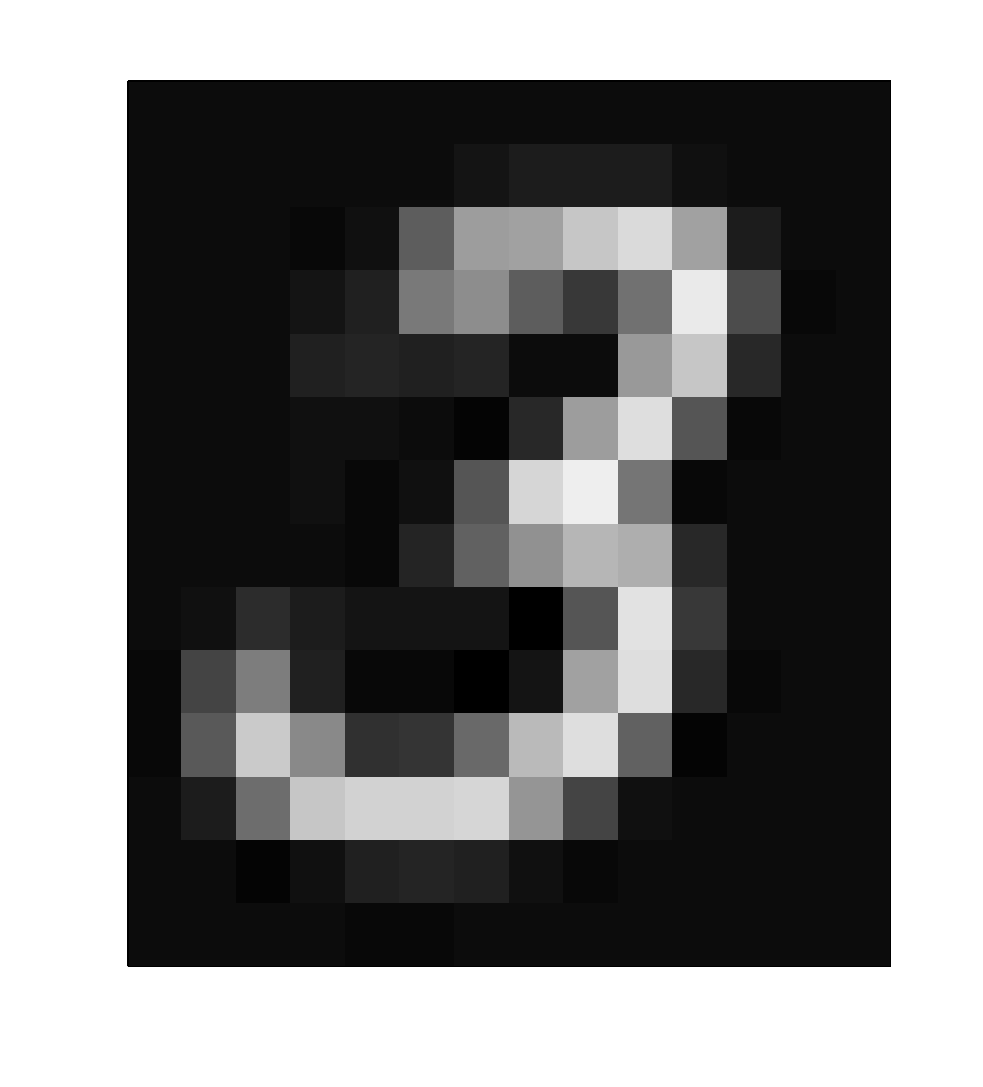}
	\end{subfigure}
	\begin{subfigure}[t]{0.89in}
		\centering
\caption*{\scriptsize Error = 0.25}
\vspace{-0.09in}
		\includegraphics[width=0.89in,
trim=0.5in 0.4in 0.35in 0.3in,clip=true]{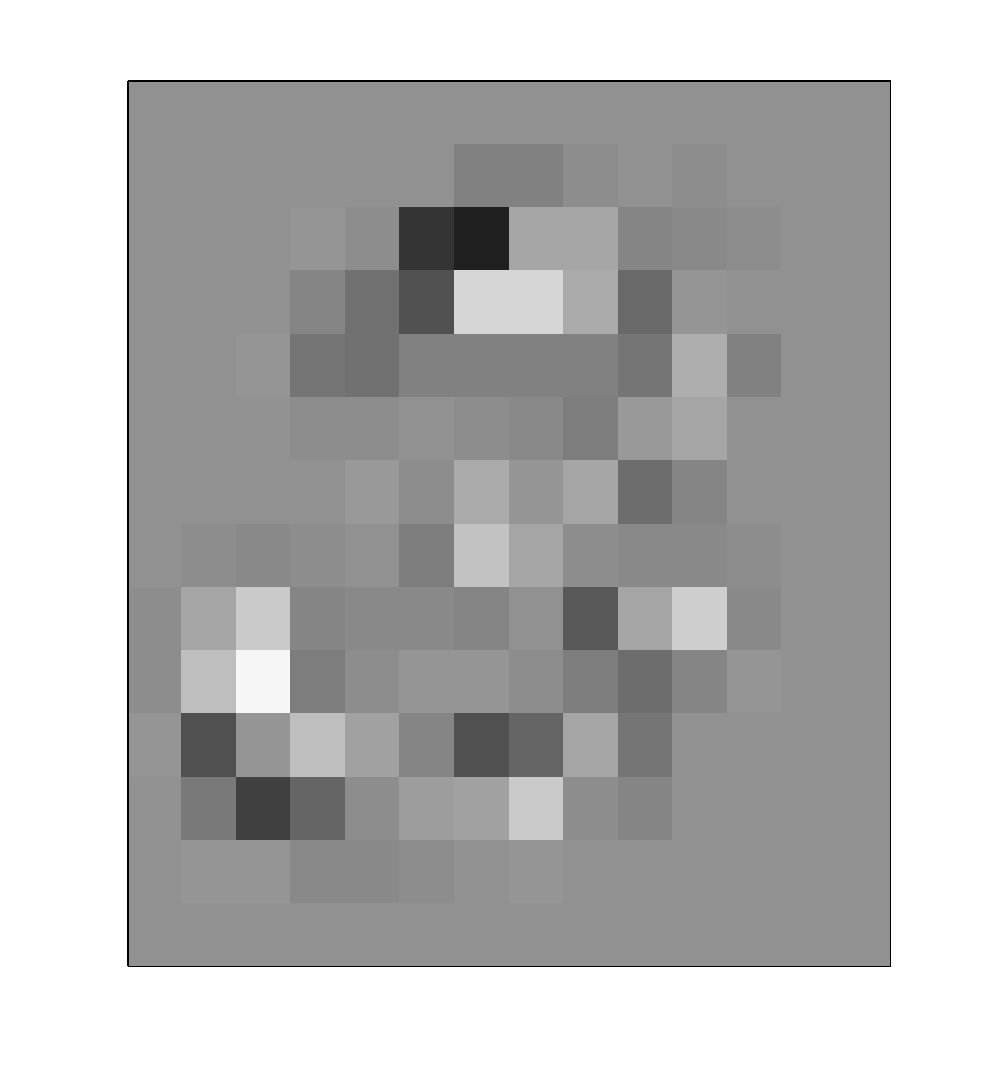}
	\end{subfigure}
	\begin{subfigure}[t]{0.89in}
		\centering
\caption*{\scriptsize Shepard \tiny ($\varepsilon = 1/\bar
h_{local}$)}\vspace{-0.082in}
		\includegraphics[width=0.89in,
trim=0.5in 0.4in 0.35in 0.3in,clip=true]{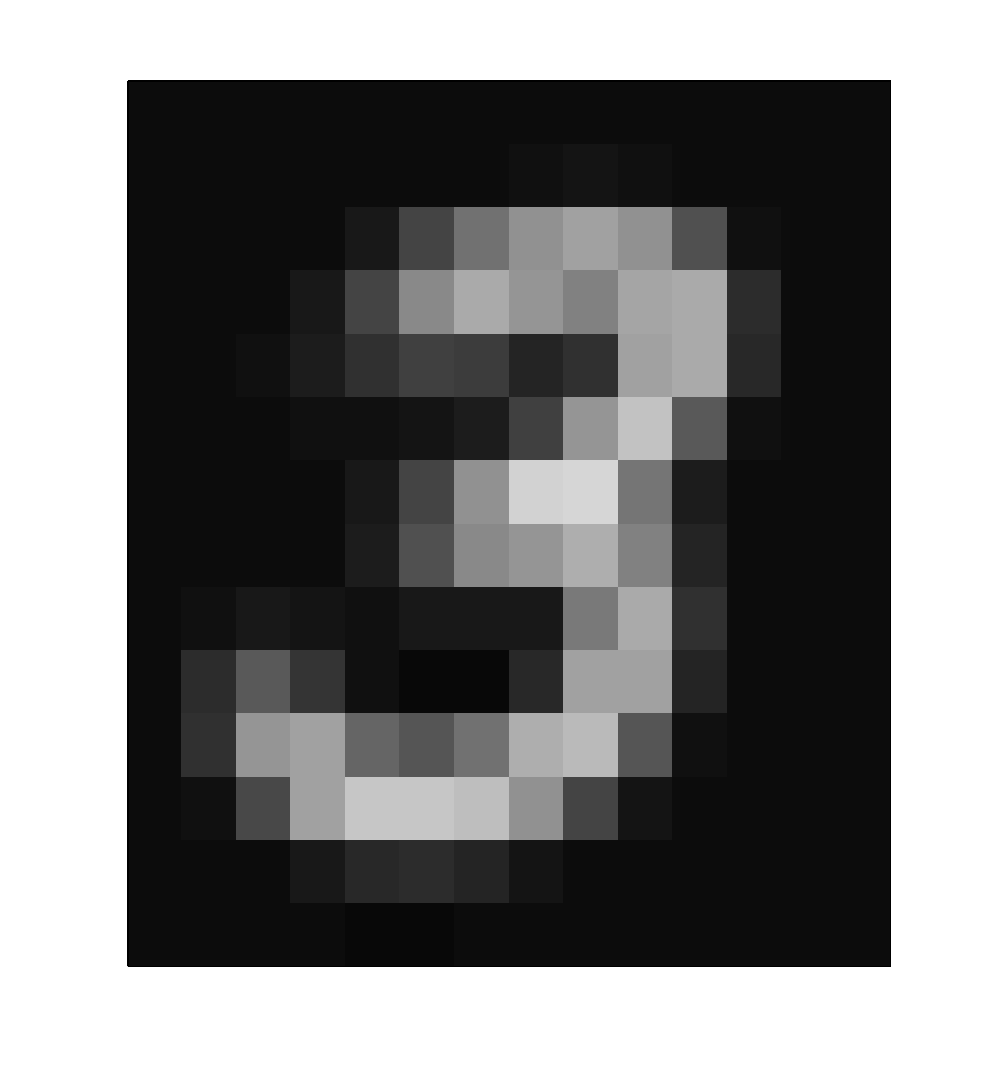}
	\end{subfigure}
	\begin{subfigure}[t]{0.89in}
		\centering
\caption*{\scriptsize Error = 0.35}\vspace{-0.09in}
		\includegraphics[width=0.89in,
trim=0.5in 0.4in 0.35in 0.3in,clip=true]{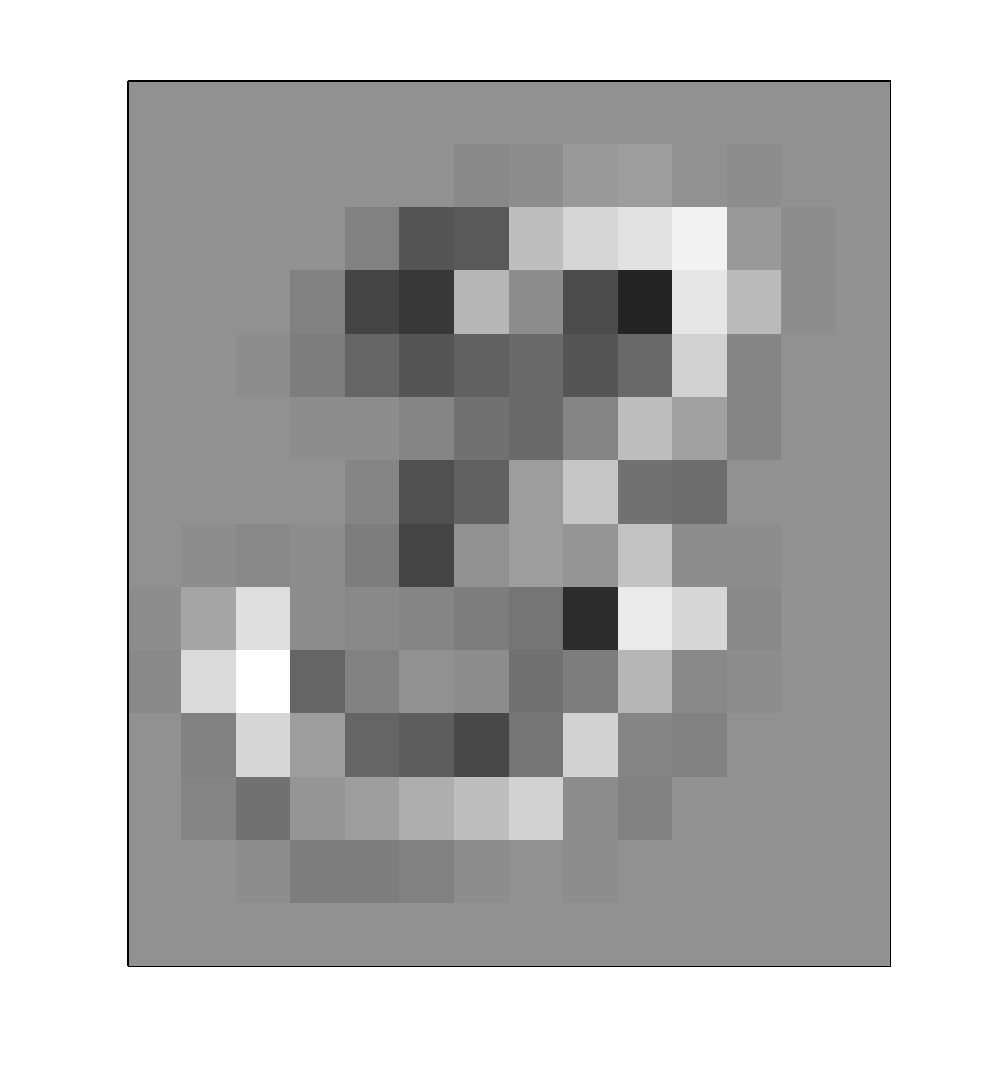}
	\end{subfigure}
	\caption{From left to right: original image to be reconstructed;
     reconstructions using the different methods, each followed by the
     residual error: cubic RBF, Gaussian RBF, and Shepard's
     method.}\label{three_example}
\end{figure}

\subsection{Handwritten Digits Datasets}
In addition to the previous synthetic example, the performance of the
inverse mapping algorithm was also assessed on a ``naturally
occurring'' high-dimensional data set: a set of digital images of
handwritten digits.  The data set (obtained from the MNIST database
\cite{digits}) consists of 1,000 handwritten images of the digits 0 to
9. The images were originally centered and normalized to have size
$28\times28$.  In our experiments, the images were further resized to
$14\times14$ pixels and normalized to have unit $l^2$ norm. We
obtained 10 different datasets, each consisting of 1,000 points in
$\R^{196}$.  The dimension reduction and subsequent leave-one-out
reconstruction were conducted on the dataset corresponding to a
specific digit, independently of the other digits. For each digit, a
10-dimensional representation of the 1,000 images was generated using
Laplacian Eigenmaps \cite{Coifman06b}.  Then the inverse mapping
techniques were evaluated on all images in the
set. Table~\ref{digits_table} shows the reconstruction error
$E_{\text{avg}}$ for the three methods, for all
digits. Fig.~\ref{three_example} shows%

\begin{figure}	
	\centering
	\begin{subfigure}[t]{0.89in}
		\centering
\caption*{\scriptsize Original}\vspace{-0.09in}
		\includegraphics[width=0.89in,
trim=0.5in 0.4in 0.35in 0.3in,clip=true]{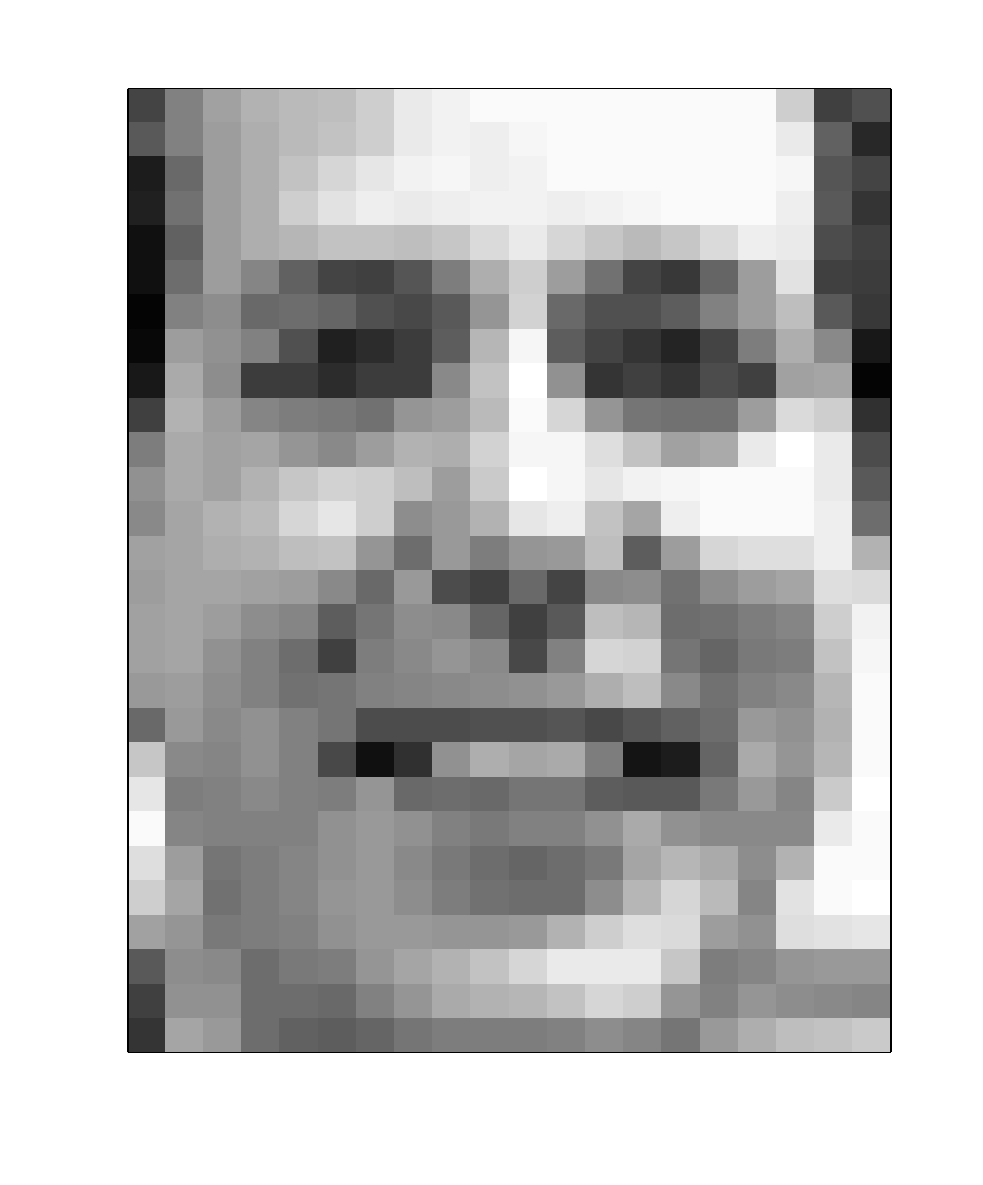}
	\end{subfigure}
	\begin{subfigure}[t]{0.89in}
		\centering
\caption*{\scriptsize Cubic RBF}\vspace{-0.09in}
		\includegraphics[width=0.89in,
trim=0.5in 0.4in 0.35in 0.3in,clip=true]{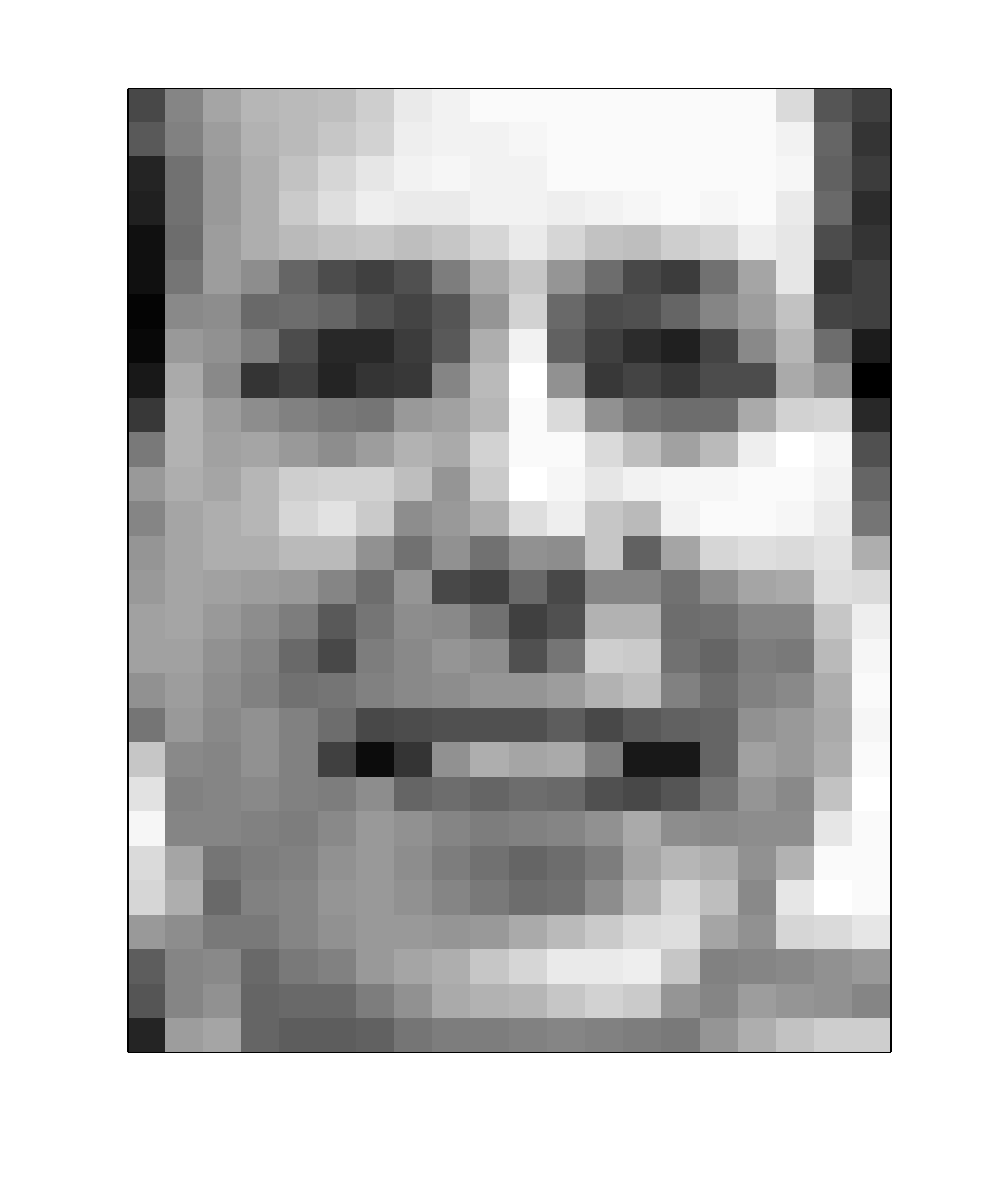}
	\end{subfigure}
	\begin{subfigure}[t]{0.89in}
		\centering
\caption*{\scriptsize Error = 0.026}\vspace{-0.09in}
		\includegraphics[width=0.89in,
trim=0.5in 0.4in 0.35in 0.3in,clip=true]{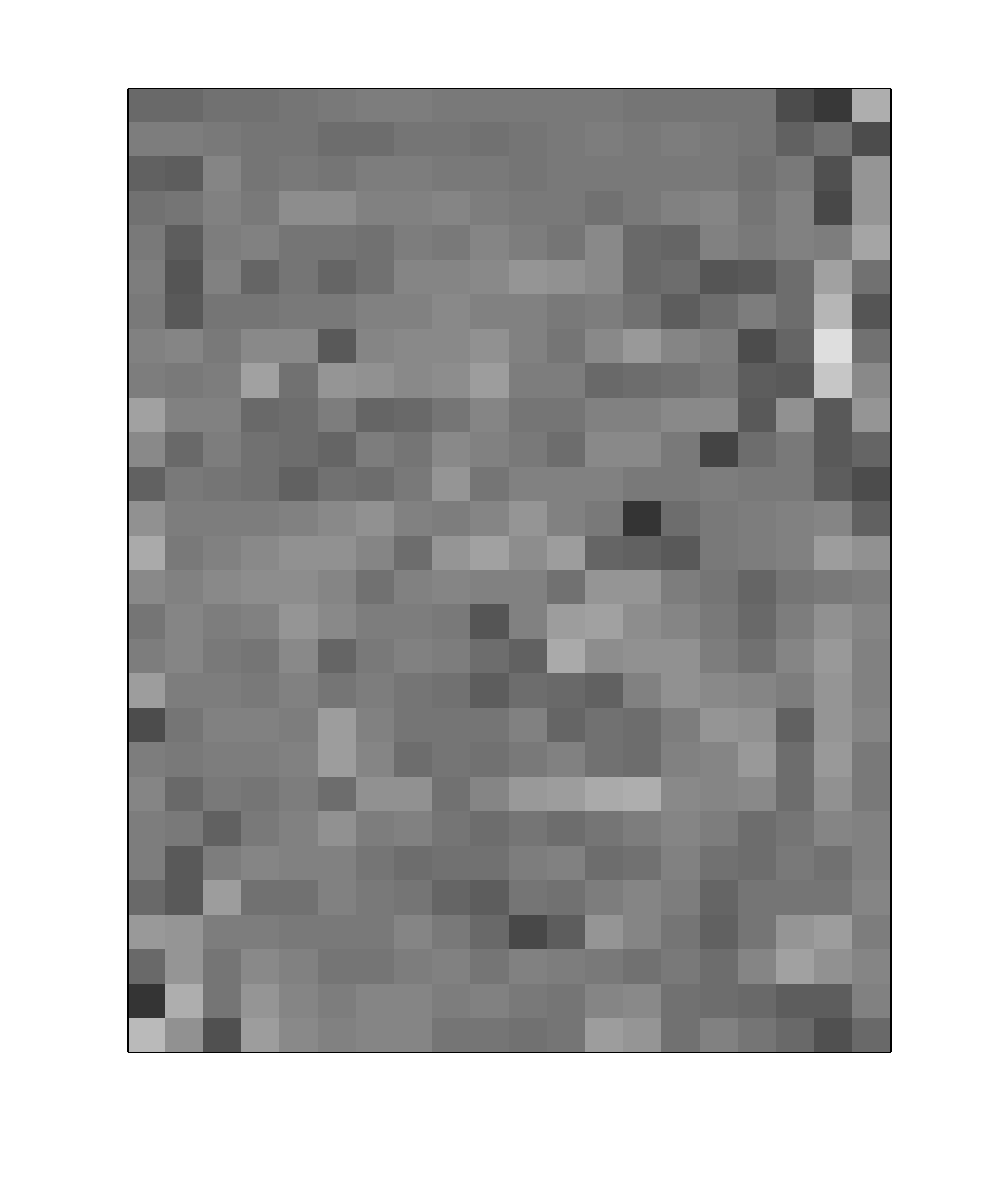}
	\end{subfigure}
	\begin{subfigure}[t]{0.89in}
		\centering
\caption*{\scriptsize Gaussian \tiny ($\varepsilon = 0.5/\bar
h_{local}$)}\vspace{-0.082in}
		\includegraphics[width=0.89in,
trim=0.5in 0.4in 0.35in 0.3in,clip=true]{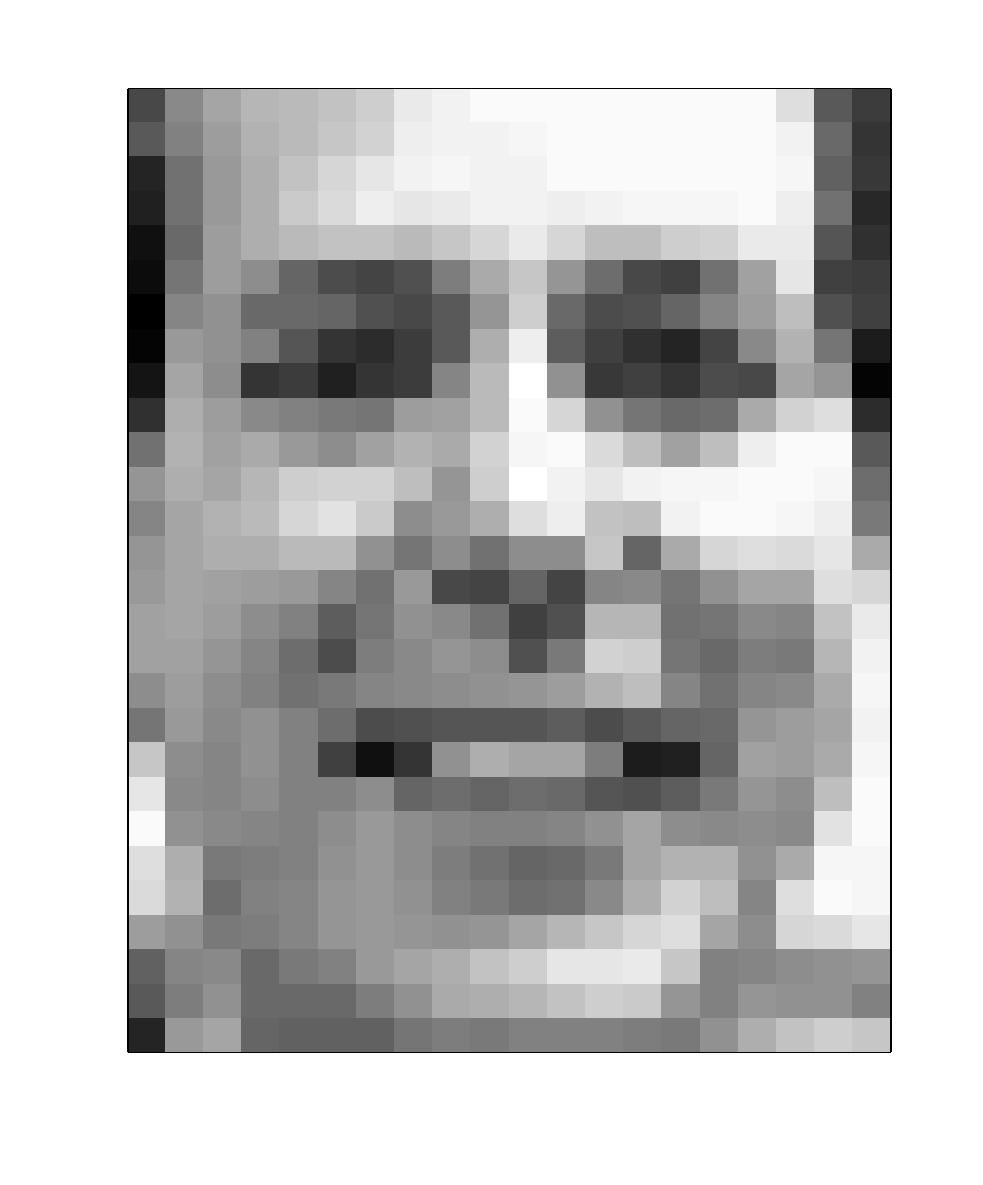}
	\end{subfigure}
	\begin{subfigure}[t]{0.89in}
		\centering
\caption*{\scriptsize Error = 0.029}
\vspace{-0.09in}
		\includegraphics[width=0.89in,
trim=0.5in 0.4in 0.35in 0.3in,clip=true]{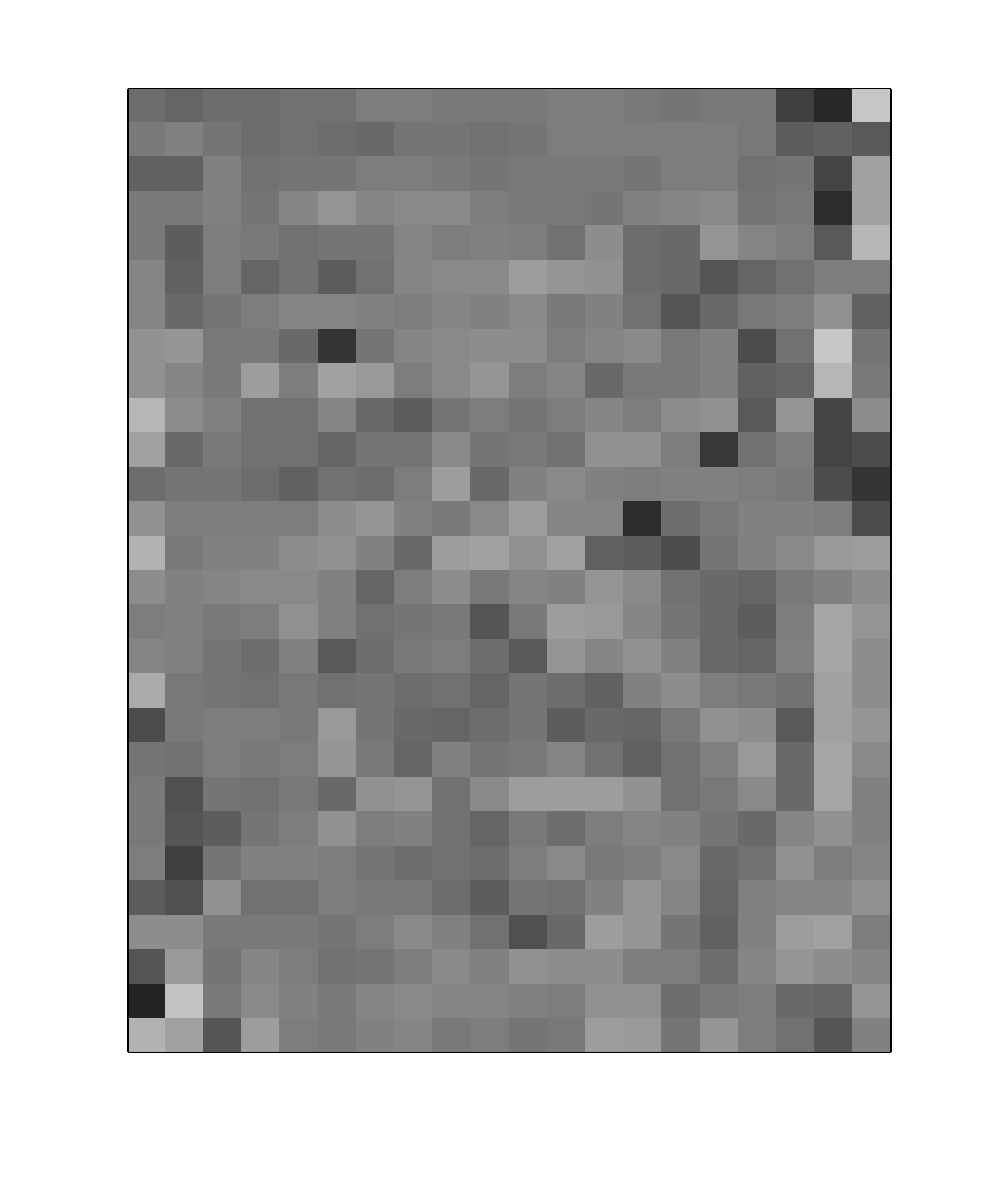}
	\end{subfigure}
	\begin{subfigure}[t]{0.89in}
		\centering
\caption*{\scriptsize Shepard \tiny ($\varepsilon = 2/\bar
h_{local}$)}\vspace{-0.082in}
		\includegraphics[width=0.89in,
trim=0.5in 0.4in 0.35in 0.3in,clip=true]{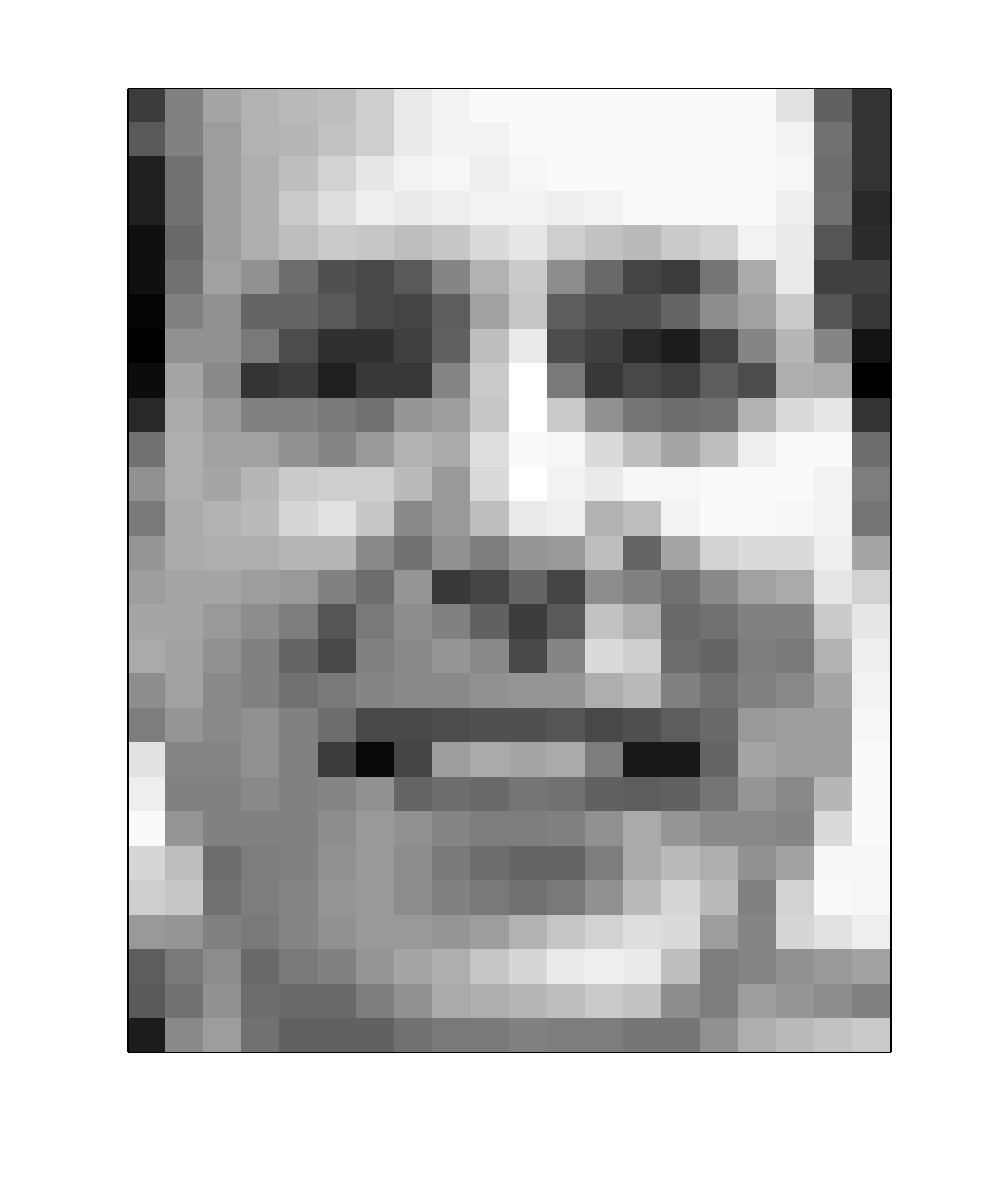}
	\end{subfigure}
	\begin{subfigure}[t]{0.89in}
		\centering
\caption*{\scriptsize Error = 0.040}\vspace{-0.09in}
		\includegraphics[width=0.89in,
trim=0.5in 0.4in 0.35in 0.3in,clip=true]{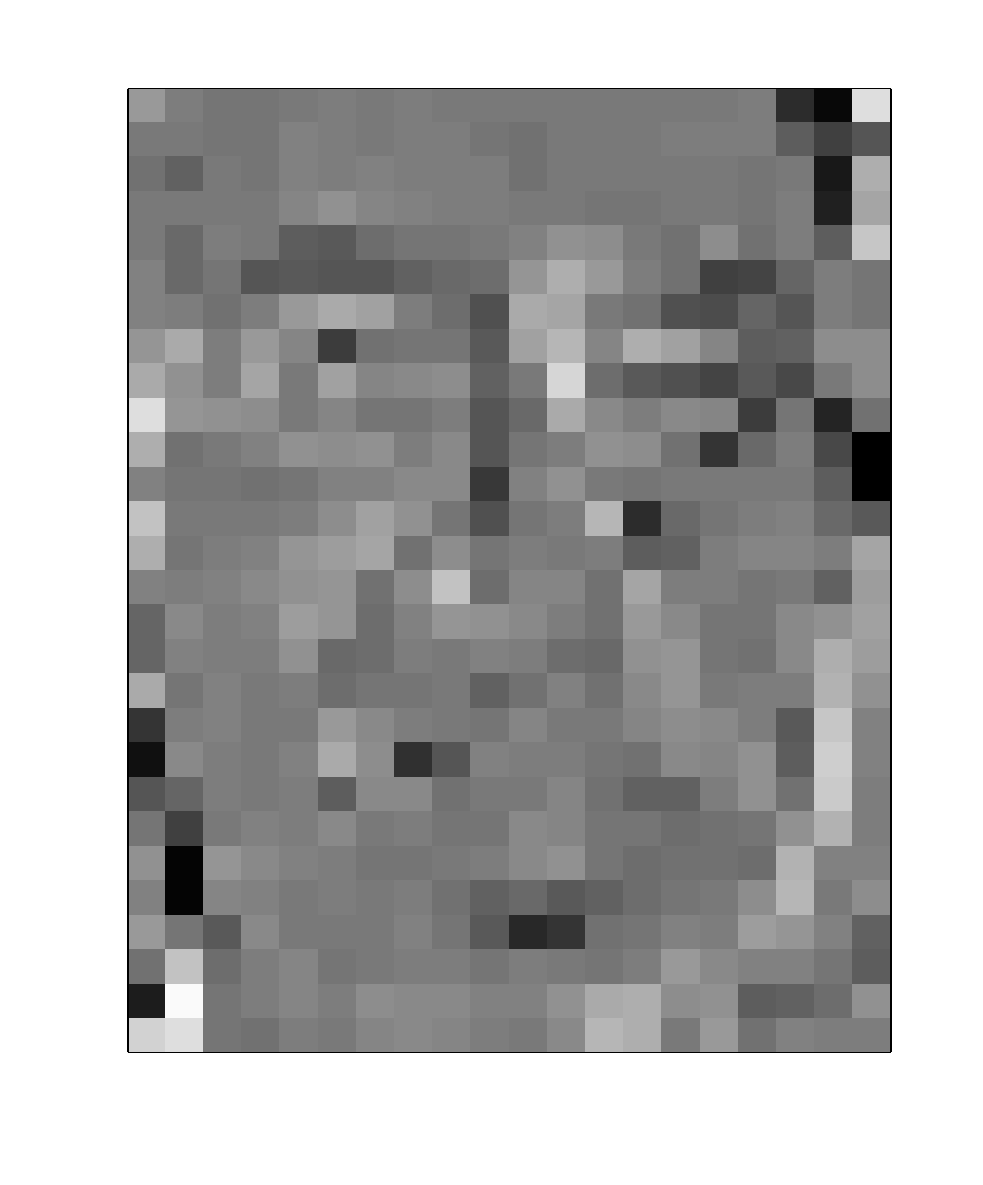}
	\end{subfigure}
	\caption{From left to right: original image to be reconstructed;
     reconstructions using the different methods, each followed by the
     residual error: cubic RBF, Gaussian RBF, and Shepard's
     method.}\label{frey_example}
\end{figure}

\begin{table}[H]
  \begin{center}
    \begin{small}
      \begin{tabular}{lllllllllll}
        \toprule
        & Cubic && \multicolumn{3}{c}{Gaussian}& & \multicolumn{3}{c}{Shepard}\\
        \midrule
        {\bf Scale ($\varepsilon \times \hloc$)}  & ---  & 
        & $0.25$ &$0.5$ & $1$ && $1$  & $2$  & $4$\\
        \midrule
        {\bf $E_{\text{avg}}$} & \red 0.0361 && 0.0457 & 0.0414  & 0.0684& & 0.0633 & 0.0603  &0.0672 \\
        \bottomrule
      \end{tabular}
    \end{small}
  \end{center}
  \caption{Reconstruction error $E_{\text{avg}}$ for the Frey Face dataset.  Red denotes lowest average
    reconstruction residual.}
  \label{frey_table}
\end{table}
\noindent  three representative
reconstructions for the digit ``3''. The optimal scales (according to
Table \ref{digits_table}) were chosen for both the Gaussian RBF and
Shepard's methods. The cubic RBF outperforms the Gaussian RBF and
Shepard's method in all cases, with the lowest average error
(Table~\ref{digits_table}), and with the most ``noise-like''
reconstruction residual (Fig.~\ref{three_example}).  Results suggest
that a poor choice of scale parameter with the Gaussian can corrupt the
reconstruction.  The scale parameter in Shepard's method must be
carefully selected to avoid the two extremes of either reconstructing
solely from a single nearest neighbor, or reconstructing a blurry, equally
weighted, average of all neighbors.
\subsection{Frey Face Dataset}
Finally, the performance of the inverse mapping algorithms was also
assessed on the Frey Face dataset \cite{freyface}, which consists of
digital images of Brendan Frey's face taken from sequential frames of
a short video. The dataset is composed of $20\times28$ gray scale
images.  Each image was normalized to have unit $l^2$ norm, providing
a dataset of 1,965 points in $\R^{560}$. A 15-dimensional
representation of the Frey Face dataset was generated via Laplacian
eigenmaps.  The inverse mapping techniques were tested on all images
in the set. Table~\ref{frey_table} shows the mean leave-one-out
reconstruction errors for the three methods. Fig.~\ref{frey_example}
shows three representative reconstructions using the different
techniques.  The optimal scales (according to
Table \ref{frey_table}) were chosen for both the Gaussian RBF and
Shepard's methods.  Again, the cubic RBF outperforms the
Gaussian RBF and Shepard's method in all cases, with the lowest
average error (Table~\ref{frey_table}), and with the most
``noise-like'' reconstruction residual (Fig.~\ref{frey_example}).
\section{Revisiting Nystr\"om}
\label{discussion}
Inspired by the RBF interpolation method, we provide in the following
a novel interpretation of the Nystr{\"o}m extension: the Nystr{\"o}m
extension interpolates the eigenvectors of the (symmetric) normalized
Laplacian matrix using a slightly modified RBF interpolation scheme.
While several authors have mentioned the apparent similarity of
Nystr{\"o}m method to RBF interpolation, the novel and detailed analysis
provided below provides a completely new insight into the limitations
and potential pitfalls of the Nystr{\"o}m extension.

Consistent with Laplacian eigenmaps, we consider the symmetric
normalized kernel $\widetilde K = D^{-1/2} K D^{-1/2}$, where $K_{ij} =
k(\xx^{(i)},\xx^{(j)})$ is a radial function measuring the similarity
between $\xx^{(i)}$ and $\xx^{(j)}$, and $D$ is the degree matrix
(diagonal matrix consisting of the row sums of $K$).  
\begin{figure}[H]
  \begin{center}
    \includegraphics[width=0.8\textwidth]{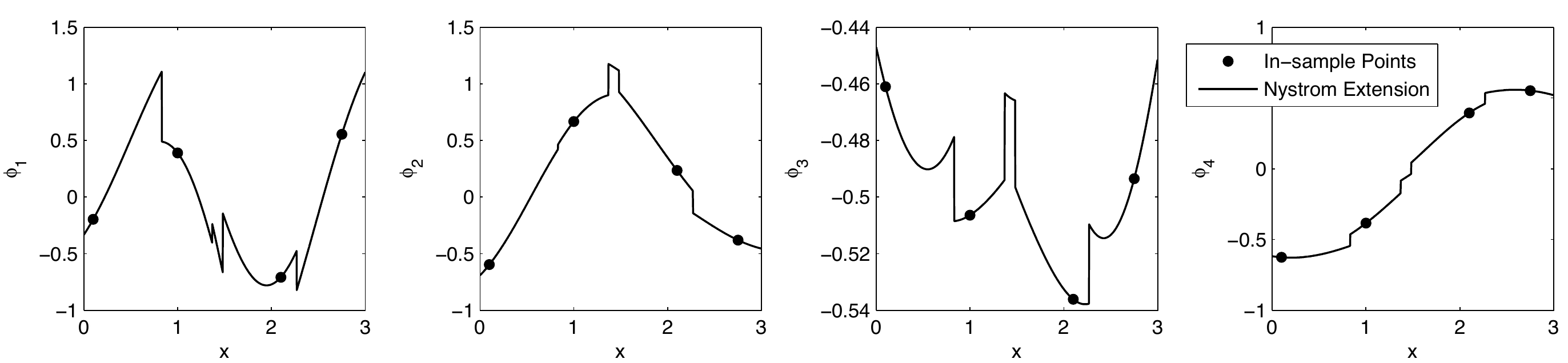}
  \end{center}
  \caption{Example of Nystr{\"o}m extension of the eigenvectors of a
    thresholded Gaussian affinities matrix.}
  \label{truncated}
\end{figure}
\noindent Given an eigenvector $\bfi_l$ of $\widetilde{K}$
(associated with a nontrivial eigenvalue $\lambda_l$) defined
over the points $\xx^{(i)}$, the Nystr{\"o}m extension of $\bfi_l$ to
an arbitrary new point $\xx$ is given by the interpolant
\begin{equation}
  \phi_l(\xx) = \frac{1}{\lambda_l} \sum_{j=1}^n \tilde k(\xx,\xx^{(j)})
  \phi_l(\xx^{(j)}),
  \label{Nystrom1}
\end{equation}
where $\phi_l(\xx^{(j)})$ is the coordinate $j$ of the eigenvector
$\bfi_l = \begin{bmatrix}\phi_l (\xx^{(1)}) & \cdots & \phi_l(\xx^{(n)})
\end{bmatrix}^T$. We now proceed by re-writing $\phi_l(\xx)$ in
(\ref{Nystrom1}), using the notation $\tilde
\kk(\xx,\cdot)=\begin{bmatrix}\tilde k(\xx,\xx^{(1)}) & \cdots &
  \tilde k(\xx,\xx^{(n)})\end{bmatrix}^T$, where $\tilde
k(\xx,\xx^{(j)}) = k(\xx,\xx^{(j)}) /\sqrt{d(\xx)d(\xx^{(j)})}$, and
$d(\xx) = \sum_{i=1}^n k(\xx,\xx^{(i)})$.
\begin{equation}
  \begin{split}
    \phi_l(\xx) & =  \lambda_l^{-1} \tilde \kk(\xx,\cdot)^T \bfi_l  =
    \tilde \kk(\xx,\cdot)^T \bFi \Lambda^{-1} [ 0 \ldots 1 \ldots 0]^T 
    =  \tilde \kk(\xx,\cdot)^T \bFi \Lambda^{-1} \bFi^T \bfi_l  \\ & =
    \tilde \kk(\xx,\cdot)^T \widetilde K^{-1} \bfi_l  =  
    \frac{1}{\sqrt{d(\xx)}} 
    \begin{bmatrix}
      k(\xx,\xx^{(1)}) & \ldots & k(\xx,\xx^{(n)}) 
    \end{bmatrix}
    D^{-1/2}  (D^{1/2} K^{-1} D^{1/2}) \bfi_l\\
    & =  \frac{1}{\sqrt{d(\xx)}}
    \kk(\xx,\cdot)^T K^{-1} ( D^{1/2} \bfi_l ) . 
  \end{split}
  \label{Nystrom2}
\end{equation}
If we compare the last line of (\ref{Nystrom2}) to (\ref{rbf_eqn}), we
conclude that in the case of Laplacian Eigenmaps, with a nonsingular kernel
similarity matrix $K$, the Nystr{\"o}m extension is computed using a
radial basis function interpolation of $\bfi_l$ after a pre-rescaling
of $\bfi_l$ by $D^{1/2}$, and post-rescaling by $1/\sqrt{d(\xx)}$.
Although the entire procedure it is not exactly an RBF interpolant, it
is very similar and this interpretation provides new insight into some
potential pitfalls of the Nystr{\"o}m method.

The first important observation concerns the sensitivity of the
interpolation to the scale parameter in the kernel $k$.  As we have
explained in section \ref{conditioning}, the choice of the optimal
scale parameter $\varepsilon$ for the Gaussian RBF is quite difficult.
In fact, this issue has recently received a lot of attention
(e.g. \cite{bermanis,coifman}).  The second observation involves the
dangers of sparsifying the similarity matrix.  In many nonlinear
dimensionality reduction applications, it is typical to sparsify the
kernel matrix $K$ by either thresholding the matrix, or keeping only the
entries associated with the nearest neighbors of each $\xx_j$.  If the
Nystr{\"o}m extension is applied to a thresholded Gaussian kernel
matrix, then the components of $\kk(\xx,\cdot)$ as well as
$\sqrt{d(\xx)}$ are discontinuous functions of $\xx$.  As a result,
$\bfi_l(\xx)$, the Nystr{\"o}m extension of the eigenvector $\bfi_l$
will also be a discontinuous function of $\xx$, as demonstrated in
Fig.~\ref{truncated}.  In the nearest neighbor approach, the extension of the
kernel function $\tilde k$ to a new point $\xx$ is highly unstable and poorly
defined.  Given this larger issue, the Nystr{\"o}m extension should not be used
in this case. In order to interpolate eigenvectors of a sparse similarity
matrix, a better interpolation scheme such as a true (non-truncated) Gaussian
RBF, or a cubic RBF interpolant could provide a better alternative to
Nystr{\"o}m.  A local implementation of the interpolation algorithm may
provide significant computational savings in certain scenarios.
\section*{Acknowledgments}
\noindent The authors would like to thank the three anonymous
reviewers for their excellent comments. NDM was supported by NSF grant DMS
0941476; BF was supported
by NSF grant DMS 0914647; FGM was partially supported by NSF grant DMS 0941476,
and DOE award DE-SCOO04096.

\bibliographystyle{model1-num-names}

\pagebreak

\end{document}